\documentclass[12pt]{amsart}
\usepackage{amsmath,amsthm,amssymb}

\def\pmatrix{\left(\begin{matrix}}
\def\endpmatrix{\end{matrix}\right)}

\def\P{{\mathbb P}}
\def\Z{{\mathbb Z}}
\def\C{{\mathbb C}}

\def\de{\delta}

\def\p{\partial}

\def\t{\theta}

\def\T{\Theta}
\def\e{\varepsilon}

\def\A{{\mathcal A}}
\def\M{{\mathcal M}}

\def\H{{\mathcal H}}

\def\D{{\mathcal D}}

\def\S{{\mathcal S}}
\def\Sn{{\mathcal S}_{\rm null}}
\def\X{{\mathcal X}}

\def\tch#1#2{{\left[\begin{matrix}{#1}\\ {#2}\end{matrix}\right]}}
\def\tt#1#2{{\t\tch#1#2}}
\def\Sp{\operatorname{Sp}(g,\Z)}
\def\Sing{\operatorname{Sing}}
\def\tn{\t_{\rm null}}

\theoremstyle{plain}
\newtheorem{thm}{Theorem}
\newtheorem{lm}[thm]{Lemma}
\newtheorem{prop}[thm]{Proposition}

\newtheorem{cor}[thm]{Corollary}

\theoremstyle{definition}

\newtheorem{rem}[thm]{Remark}

\title{Singularities of the theta divisor at points of order two}
\author{Samuel Grushevsky\and Riccardo Salvati Manni}
\address{Mathematics Department, Princeton University, Fine Hall,
Washington Road, Princeton, NJ 08544, USA. Research is supported in
part by National Science Foundation under the grant DMS-05-55867.}
\email{sam@math.princeton.edu}
\address{Dipartimento di Matematica, Universit\`a ``La Sapienza'',
Piazzale A. Moro 2, Roma, I 00185, Italy}
\email{salvati@mat.uniroma1.it}
\date{\today}
\begin{document}
\begin{abstract}
In this note we study the geometry of principally polarized abelian
varieties (ppavs) with a vanishing theta-null (i.e. with a singular
point of order two and even multiplicity lying on the theta divisor)
--- denote $\tn$ the locus of such ppavs. We describe the locus
$\tn^{g-1}\subset\tn$ where this singularity is not an ordinary
double point. By using theta function methods we first show
$\tn^{g-1}\subsetneq\tn$ (this was shown in
\cite{debarredecomposes}, see below for a discussion). We then show
that $\tn^{g-1}$ is contained in the intersection $\tn\cap N_0'$ of
the two irreducible components of the Andreotti-Mayer
$N_0=\tn+2N_0'$, and describe by using the geometry of the universal
scheme of singularities of the theta divisor which components of
this intersection are in $\tn^{g-1}$.

Some of the intermediate results obtained along the way of our proof
were concurrently obtained independently by C. Ciliberto and G. van
der Geer in \cite{amsp} and by R. de Jong in \cite{dj}, version 2.
\end{abstract}
\maketitle

\section*{Introduction}
The theta divisor $\T$ of a generic principally polarized abelian
variety (ppav) is smooth. In \cite{am} and \cite{gen4am} it was
shown that ppavs $(A,\T)$ with a singular theta divisor form a
divisor $N_0$ in the moduli space of ppavs $\A_g$.

In \cite{mu} and \cite{debarredecomposes} it was shown that $N_0$
has two irreducible components $\tn$ and $N_0'$, where $\tn$ is the
locus of those ppavs for which the theta divisor has a singularity
at a point of order two, and $N_0'$ is the closure of the locus of
those ppavs for which the theta divisor has a singularity at a point
not of order two.

Moreover it has been proven that for a generic ppav $(A,\T)\in\tn$
the theta divisor $\T$ has a unique singular point, which is a
double point. Similarly the generic element of $N_0'$ is a
principally polarized abelian variety for which the theta divisor
$\T$ has two distinct singular points $x$ and $-x$, which are double
points --- this is due to the symmetry of the theta divisor. So we
can write
\begin{equation}\label{divisor}
 N_0=\tn+2N_0 '.
\end{equation}

Moreover, in \cite{debarredecomposes} the author claims that in the
case of $\tn$, the singular point is generically an ordinary double
point (i.e. the tangent cone to the theta divisor at such a point
has maximal rank, i.e. rank $g$). He refers to \cite{te} for a
proof. It seems that the reference is not really appropriate, since
\cite{te} treats the restriction of $\tn$ to $\M_g$, the moduli
space of curves of genus $g$. But in this case, for $g\geq 4$ we
know that the double points cannot be ordinary.

Upon reading a preliminary version of this note, O. Debarre has
explained to us that this can be fixed by a little more work, using
the results in his work --- see remark \ref{debarreproof} for that
proof.  The different method we use, however, yields a further
insight into the solution of some interesting problems about the
double points on the theta divisor, and thus is hopefully of
independent interest.

\smallskip
Our interest in the divisor $\tn$ comes from our recent paper
\cite{genus4} where we gave, in genus 4, a characterization of the
intersection
\begin{equation}\label {cap}
 \tn\cap N_0'
\end{equation}
in terms of the rank of the tangent cone of the singular point of the
theta divisor. Essentially we proved that the above locus is
characterized by the fact that the singular point of order two is a
double point, but not an ordinary double point. The existence of
ppavs of dimension 4 for which the theta divisor has an ordinary
double point was proved in \cite{gen4am}.

We also study the intersection (\ref{cap}) and its relation with the
locus $\tn^{g-1}$ --- the sublocus of $\tn$ parametrizing ppavs $(A,
\T)$ with a singular point of order two which is a double point, but
not an ordinary double point.

\smallskip
Incidentally with our proof, we solve a question raised in
\cite{dj}, about the vanishing locus of a modular form defined on
the universal theta divisor, which is constructed in that work. We
were recently informed by R. de Jong that he has also obtained
independently a proof of proposition \ref{gauss} in the new version
of his work.

\smallskip
At the last stages of editing the draft of this text we were
informed by C. Ciliberto and G. van der Geer of their preprint
\cite{amsp} where along the way of their discussion of the dimension
of Andreotti-Mayer loci they also discuss related questions about
the singularities of theta divisors for points on $\tn$. Their main
interest is in the higher Andreotti-Mayer loci, and our results on
$\tn^{g-1}$ do not have a parallel in their work.

\section{Notations and definitions}
In this section we recall notations, definitions, as well as some
results from \cite{genus4}. We denote $\H_g$ the {\it Siegel upper
half-space}, i.e. the set of symmetric complex $g\times g$ matrices
$\tau$ with positive definite imaginary part. Each such $\tau$
defines a complex principally polarized abelian variety (ppav for
short) $\C^g/\tau\Z^g+\Z^g$. If $\sigma=\pmatrix a&b\\
c&d\endpmatrix\in\Sp$ is a symplectic matrix in a $g\times g$ block
form, then its action on $\tau\in\H_g$ is defined by
$\sigma\cdot\tau:=(a\tau+b)(c\tau+d)^{-1}$, and the moduli space of
ppavs is the quotient $\A_g=\H_g/\Sp$. A period matrix $\tau$ is
called {\it decomposable} if there exists $\sigma\in\Sp$ such that
$$
 \sigma\cdot\tau=\pmatrix \tau_1&0\\
 0&\tau_2\endpmatrix,\quad\tau_i\in\H_{g_i},\ g_1+g_2=g, g_i>0;
$$
otherwise we say that $\tau$ is indecomposable.

For $\e,\de\in (\Z/2\Z)^g$, thought of as vectors of zeros and ones,
$\tau\in\H_g$ and $z\in \C^g$, the {\it theta function with
characteristic $[\e,\de]$} is
$$
 \tt\e\de(\tau,z):=\sum\limits_{m\in\Z^g} \exp \pi i \left[
 ^t(m+\frac{\e}{2})\tau(m+\frac{\e}{2})+2\ ^t(m+\frac{\e}{2})(z+
 \frac{\de}{2})\right].
$$
Sometimes we shall write $\t(\tau, z)$ for the theta function with
characteristic $[0,0]$.

A {\it characteristic} $[\e,\de]$ is called {\it even} or {\it odd}
depending on whether $\tt\e\de(\tau,z)$ is even or odd as a function
of $z$, which corresponds to the scalar product
$\e\cdot\de\in\Z/2\Z$ being zero or one, respectively. A {\it theta
constant} is the evaluation at $z=0$ of a theta function. All odd
theta constants of course vanish identically in $\tau$.

Observe that
$$
 \tt{0}{0}\left(\tau,z+\tau\frac\e2+\frac\de2\right)= \exp \pi i
 \left(-\frac{^t\e}{2}\,\tau\frac{\e}{2}\,\,-\frac{^t\e}{2}\,(z+
 \frac{\de}{2})\right)\tt\e\de(\tau,z),
$$
i.e. theta functions with characteristics are, up to some non-zero
factor, the Riemann theta function (the one with characteristic
$[0,0]$) shifted by points of order two.

A map $f:\H_g\to \C$ is called a {\it (scalar) modular form} of
weight $k$ with respect to a finite index
subgroup $\Gamma\subset\Sp$ if
$$
 f(\sigma\cdot\tau)={\rm det}(c\tau+d)^k f(\tau)\qquad\forall
 \tau\in\H_g,\forall\sigma\in\Gamma,
$$
and if additionally $f$ is holomorphic at all cusps of
$\H_g/\Gamma$.

We define
$$
 \Gamma_g(n):=\left\lbrace \sigma=\pmatrix a&b\\ c&d\endpmatrix
 \in\Sp \, |\, \sigma\equiv\pmatrix 1&0\\
 0&1\endpmatrix\ {\rm mod}\ n\right\rbrace
$$
$$
 \Gamma_g(n,2n):=\left\lbrace \sigma\in\Gamma_g(n)\, |\, {\rm
 diag}(a^tb)\equiv{\rm diag} (c^td)\equiv0\ {\rm mod}\
 2n\right\rbrace.
$$
These are finite index normal subgroups of $\Sp$.

Under the action of $\sigma\in\Sp$ the theta functions transform as
follows:
$$
 \t\bmatrix \sigma\pmatrix \e\\ \de\endpmatrix\endbmatrix
 (\sigma\cdot\tau,\,^{t}(c\tau+d)^{-1}z)\qquad\qquad\qquad
$$
$$
 \qquad\qquad\qquad=\phi(\e,\,\de,\,\sigma,\,
 \tau,\,z)\det(c\tau+d)^{\frac{1}{2}}\t\bmatrix \e\cr \de
 \endbmatrix(\tau,\,z),
$$
where
$$
 \sigma\pmatrix \e\cr \de\endpmatrix :=\pmatrix d&-c\cr
 -b&a\endpmatrix\pmatrix \e\cr \de\endpmatrix+ \pmatrix {\rm
 diag}(c \,^t d)\cr {\rm diag}(a\,^t b)\endpmatrix,
$$
considered in $(\Z/2\Z)^g$, and $\phi(\e,\,\de,\,\sigma,\,\tau,\,z)$
is some complicated explicit function. For more details, we refer to
\cite{igbook} and \cite{farkasrauch}.

Thus theta constants with characteristics are modular forms of
weight $1/2$ with respect to $\Gamma_g(4,8)$, i.e. we have
$$
 \tt\e\de(\sigma\cdot\tau,0)=\det(c\tau+d)^{1/2}\tt\e\de(\tau,0)
 \qquad \forall \sigma\in\Gamma_g(4,8).
$$

The theta constants are known to define an embedding of the {\it
level moduli space} $\A_g(4,8):=\H_g/\Gamma_g(4,8)$ (see
\cite{igbook}, chapter V):
$$
 \begin{aligned}
 Th:\A_g(4,8)&\rightarrow\P^{2^{g-1}(2^g+1)-1}&\\
 \tau&\mapsto\left\lbrace\tt\e\de(\tau)\right\rbrace_{[\e,\de]\ {\rm
 even}}& \end{aligned},
$$
which extends to the Satake compactification
${\overline{\A_g(4,8)}}$.

We call the {\it theta-null divisor} $\tn\subset\A_g$ the zero locus
of the product of all even theta constants. We define a
stratification of $\tn$ as follows. For $h=0,\ldots,g$ we let
$$
 \tn^h=\left\lbrace\tau\in\H_g: \exists[\e,\de]\ {\rm even}\ ,
 \tt\e\de(\tau)=0;\ \operatorname{rk}\left(
 \frac{\partial^2\tt\e\de(\tau,z)}{\partial z_j\partial
 z_k}\right)_{z=0}\le h\right\rbrace,
$$
i.e. the locus of points on $\tn$ where the rank of the tangent cone
to the theta divisor at the corresponding point
$\frac{\tau\e+\de}{2}$ of order two is at most $h$.

The partial derivatives of the theta function are not modular forms.
However, since theta functions satisfy the heat equation
$$
 \frac{\partial^2\tt\e\de(\tau,z)}{\partial z_j\partial z_k}
 =\pi
 i(1+\delta_{j,k})\frac{\partial\tt\e\de(\tau,z)}{\partial\tau_{jk}},
$$
(where $\delta_{j,k}$ is Kronecker's symbol) and the partial
$\tau_{jk}$ derivative of a section $\tt\e\de(\tau,0)$ of a line
bundle on $\A_g(4,8)$, when restricted to the zero locus of this
theta constant, is a section of the same bundle, on the locus
$\lbrace\tt\e\de(\tau,0)=0\rbrace$ the second derivative
$$
 \frac{\partial^2\tt\e\de(\tau,z)}{\partial z_j\partial z_k}|_{z=0}
$$
is a modular form for $\Gamma_g(4,8)$.

Using the heat equation, the Hessian of the theta function with
respect to $z$ can be rewritten using the first derivatives with
respect to $\tau_{jk}$. Hence if a point $x=\tau\frac\e2+\frac\de2$
of order two is a singular point on the theta divisor, which is
simply to say $\tt 00(\tau,x)=0=\tt\e\de(\tau,0)$ (the first
derivatives at zero of an even function are all zero), the rank of
the quadric defining the tangent cone at $x$ is the rank of the
matrix obtained by applying the $g\times g$-matrix-valued
differential operator
$$
 \D:= \left(\begin{array}{rrrr}
 \,\frac{\p}{\p\tau_{11}}&\frac{1}{2}\frac{\p}
 {\p\tau_{12}}&\dots&\frac{1}{2}\frac{\p}{\p\tau_{1 g}}\\
 \frac{1}{2}\frac{\p}{\p \tau_{21}}&\frac{\p}{\p
 \tau_{22}}&\dots&\frac{1}{2}\frac{\p}{\p\tau_{2 g}}\\
 \dots&\dots&\dots&\dots\\
 \frac{1}{2}\frac{\p}{\p \tau_{g 1}}& \dots&\dots& \,\
 \frac{\p}{\p\tau_{g g}}\end{array}\right)
$$
to $\tt\e\de(\tau,0)$.

In $\A_g(4,8)$, the locus $\tn^{h}$ is given by the conditions
\begin{equation}\label {char}
 \lbrace \exists\,\, [\e,\de]\ {\rm even;}\ 0=\tt\e\de
 (\tau);\ {\rm rk}\, \D\tt\e\de(\tau)\le h\rbrace.
\end{equation}
The divisor $\tn\subset \A_g(4,8)$ is reducible. Its irreducible
components are the divisors of individual theta constants with
characteristics (cf. \cite{fr} page 88 for $g\geq 3$; it is easily
verified also for $g=1,2$). These components are all conjugate under
the action of $\Sp/\Gamma(4,8)$, and it follows that $\tn$ and
$\tn^{h}$ are well-defined on $\A_g$ and not only on $\A_g(4,8)$.

\section{Partial toroidal compactification}
In this section we recall from \cite{mu} the partial toroidal
compactification of the moduli space of ppavs, and the description
of the intersection of a subvariety of $\A_g$ with the boundary
$\partial\A_g$. We will be mainly interested in the intersection of
the divisor $\tn$ with the boundary.

The partial compactification that we consider is
\begin{equation}\label{tor}
 \overline{\A_g}^1:=\A_g\cup\partial\A_g.
\end{equation}
This is the coarse moduli space of ppavs $(A,\T)$ of dimension $g$
and their rank 1 degenerations, obtained as the blowup of
$\A_g\sqcup\A_{g-1}$ along $\A_{g-1}$.

We denote $p:\overline{\A_g}^1\to\A_g\sqcup\A_{g-1}$ the projection
map. An element of $\p\A_g$ is a pair $(\overline{G}, D)$, where
$\overline{G}$ is a complete $g$-dimensional variety that is a limit
of $g$-dimensional abelian varieties, and $D$ is an ample divisor
that is the limit of the respective theta divisors. Obviously an
element of $\A_{g-1}$ is a pair $(B,\Xi)$ where $B$ is a ppav of
dimension $g-1$ and $\Xi$ is its theta divisor.

The restriction of the map
\begin{equation} \label{fib}
 p_{|\partial\A_g}:\partial\A_g\to\A_{g-1}
\end{equation}
has $B/Aut (B, \Xi)$ as fiber over $(B,\Xi)$. This means that the
fiber over a general $(B,\Xi)\in\A_{g-1}$ is the Kummer variety
$B/\pm 1$.

We know from \cite{mu}
\begin{thm}\label{muthm}
\begin{equation}\label {inter}
 \tn\cap\partial\A_g=\left(\bigcup_{ (\overline{G},
 D)}2_B(\Xi)\right)\cup p^{-1}(\t_{null,\,g-1}).
 \end{equation}
where we denoted $2_B(\Xi):=\left\lbrace 2x | x\in\Xi\right\rbrace$.
\end{thm}

\section{Double points on general elements of $\tn$ are ordinary}
As we stated in the introduction we shall prove that a generic ppav
in $\tn$ has an ordinary double point as the only singular point of
the theta divisor. We need a technical lemma
\begin{lm}\label{mainlm}
Let $F(x_1, \dots, x_n)=0$ be the equation of  a hypersurface
$X\subset \C^n$. Let $F_1, \dots, F_n$ be the partial derivatives of
$F$, and let
$$
  G:\C^n\to \C^n
$$
be the  induced map
$$
 G(x_1,\dots,x_n):=\left(F_1(x_1, \dots ,x_n), \dots, F_n(x_1, \dots
 ,x_n)\right).
$$
Then $dF$ ramifies at $x\in X$ if and only if the matrix
\begin{equation}\label{hess}
 \left(\begin{array}{rrrrr}
 \,\frac{\p^2 F(x)}{\p X_1\p X_1}&\frac{\p^2 F(x)}{\p
 X_1\p X_2}
 &\dots&\frac{\p^2 F(x)}{\p X_1\p X_{n}}&\frac{\p F(x)}{\p X_1}\\
 \,\frac{\p^2 F(x)}{\p X_1\p X_2}&\frac{\p^2 F(x)}{\p
 X_2\p X_2}
 &\dots&\frac{\p^2 F(x)}{\p X_2\p X_{n}}&\frac{\p F(x)}{\p X_2}\\
  \dots&\dots&\dots&\dots&\dots\\
  \,\frac{\p^2 F(x)}{\p X_n\p X_1}&\frac{\p^2 F(x)}{\p
 X_n\p X_2}
 &\dots&\frac{\p^2 F(x)}{\p X_n\p X_{n}}&\frac{\p F(x)}{\p X_n}\\
\frac{\p F(x)}{\p X_1}& \dots&\dots&\frac{\p F(x)}{\p X_n}&0\end{array}\right)
\end{equation}
does not have maximal rank.
\end{lm}
\begin{proof}
The map $G$ ramifies  if and only if there exists a vector $v\in
T_x(X)$ of the tangent space mapping to $\lambda dF(x)$, where $dF$
is the column vectors whose  entries are the $F_i$. Denoting by $H$
the hessian matrix of $F$, this becomes
$$
  H(x)v=\lambda dF(x).
$$
Obviously, since $v\in T_x(X)$, we have the scalar product $v\cdot
dF(x)=0$. Thus the two assertions are equivalent to
$$
  \left(\begin{array}{rr}H(x)&dF(x)\\
  ^t dF(x)&0\end{array}\right)  \left(\begin{array}{r} v\\
  -\lambda \end{array}\right) =  \left(\begin{array}{r}0\\
0 \end{array}\right).
$$
But this is true if and only if  the matrix (\ref{hess}) does not
have maximal rank.
\end{proof}
Now we are able to prove the following
\begin{thm}\label{mainthm}
$$
 \tn^{g-1} \subsetneq\tn
$$
\end{thm}
\begin{proof} Let us consider the intersection of (the closure in
$\overline{\A_g}^1$ of) these loci with the boundary $\p\A_g$. We
shall restrict ourselves to considering in (\ref{inter}) the
component
$$
 \left(\bigcup_{ (\overline{G}, D)}2_B(\Xi)\right)
$$
of $\tn\cap\p\A_g$. Let $\t(\tau ', z)$ be the standard theta
function of genus $g-1$. By using the heat equation we can easily
check that $\tn^h\cap\p\A_g$ restricted to the above component is
described by the following analytic conditions
\begin{equation}
 \t(\tau', z/2)=0
\end{equation}
\begin{equation}\label{matr}
 \operatorname{rk}\left(\begin{array}{rrrrr}
 \,\frac{\p^2\t(\tau', z/2)}{\p z_1\p z_1}&\frac{\p^2\t(\tau', z/2)}{\p
 z_1\p z_2}
 &\dots&\frac{\p^2\t(\tau', z/2)}{\p z_1\p z_{g-1}}&\frac{\p\t(\tau',
 z/2)}{\p z_1}\\
 \frac{\p^2\t(\tau', z/2)}{\p z_2\p z_1}&\frac{\p^2\t(\tau', z/2)}{\p
 z_2\p z_2}
 &\dots&\frac{\p^2\t(\tau', z/2)}{\p z_2\p z_{g-1}}&\frac{\p\t(\tau',
 z/2)}{\p z_2}\\
 \dots&\dots&\dots&\dots&\dots\\
 \frac{\p^2\t(\tau', z/2)}{\p z_{g-1}\p z_1}&\frac{\p^2\t(\tau', z/2)}{\p
 z_{g-1}\p z_2}
 &\dots&\frac{\p^2\t(\tau', z/2)}{\p z_{g-1}\p
 z_{g-1}}&\frac{\p\t(\tau', z/2)}{\p z_{g-1}}\\
 \frac{\p\t(\tau', z/2)}{\p z_1}& \dots&\dots& \frac{\p\t(\tau', z/2)}{\p
 z_{g-1}}&0\end{array}\right)\le h.
\end{equation}
The fact that the  determinant of (\ref{matr}) does not vanish
identically follows immediately from the previous lemma, since it is
a well known fact that the Gauss map of the theta divisor of an
abelian variety does not ramify everywhere. But this is what we need
to finish the proof of the theorem. Indeed, if the determinant of
the above matrix is non-zero, it means that its rank is equal to
$g$, and thus the corresponding boundary point lies in
$(\tn\setminus\tn^{g-1})\cap\p\A_g$.

We now note that in \cite{dj} the determinant of matrix (\ref{matr})
restricted to $\t(\tau',z)=0$ was studied. It is denoted $\eta
(\tau',z)$ there (though explained in a slightly different way ---
see the remark below for a discussion), and it is shown there in
particular that it does not vanish identically for $\tau'$ in the
Jacobian locus.
\end{proof}
\begin{rem}
From our approach to the proof of theorem \ref{mainthm} we see the
geometric significance of the matrix (\ref{matr}), which is of
independent interest, and allows us to answer a question posed in
\cite{dj} about the vanishing of the modular form $\eta(\tau',z)$
defined on the theta divisor. Indeed, let us write the matrix in
(\ref{matr}) as
$$
 B=\left(\begin{array}{rr}
 H&dF\\
 ^t dF&0\end{array}\right),
$$
We denote by $H^c$
the matrix of cofactors of $H$. Then in \cite{dj} the definition is
$$
 \eta(\tau',z):={}^t dF H^c dF\, (\tau', z)
$$
A simple computation of the determinant of the matrix $B$, expanded
using the last line and the last column, shows that
$$
 \det B (\tau ', z)={}^t dFH^c dF (\tau ', z)
$$
\end{rem}
As an immediate  consequence of the proof of the previous theorem,
we solve a problem raised in \cite{dj} (R. de Jong has informed us
that he has also obtained a proof of this independently,
included in the updated version of \cite{dj}).
\begin{prop}\label{gauss}
The function $\eta(\tau',z)$ vanishes at the  point $(\tau_0, x_0)$
if and only if  $x_0$  is a ramification  point for the Gauss map
$G_{\tau_0}$ of the theta divisor of the abelian variety with period
matrix $\tau_0$.
\end{prop}
In \cite{dj} it is proved that $\eta(\tau',z)$ is a theta function
of order $g$ with respect to $z$ and weight $(g+4)/2$ with respect
to $\tau$ (note that in our notations there is a shift of $g$ by
$-1$ compared to \cite{dj}). In \cite{genus4} we proved that the
function
\begin{equation}\label{matr2}
F(\tau)= \operatorname{det}\left(\begin{array}{rrrr}
 \,\frac{\p\t(\tau, 0)}{\p \tau_{11}}&( \frac{\p \t(\tau, 0)}{2\p
 \tau_{1 2}}
 &\dots& \frac{\p \t(\tau, 0)}{2\p \tau_{1 g}}\\
 \frac{\p \t(\tau, 0)}{2\p \tau_{1 2}}& \frac{\p \t(\tau, 0)}{\p
 \tau_{2 2}}
 &\dots& \frac{\p \t(\tau, 0)}{2\p \tau_{2 g}}\\
\dots&\dots&\dots&\dots\\
 \frac{\p \t(\tau, 0)}{2\p \tau_{1 g}}&\dots&\dots& \frac{\p \t(\tau,
 0)}{\p \tau_{g g}}\
 \end{array}\right)
\end{equation}
defines a modular form of weight $(g+4)/2$ with respect to
$\Gamma_g(4, 8) $ along $\t(\tau, 0)=0$. Obviously the determinant
of the matrix (\ref{matr}) gives the intersection of $F(\tau)$ with
a suitable boundary component of  the partial toroidal
compactification of $A_g(4,8)$. \smallskip

We end this section by observing that in the  above discussion the
intersection $\tn^h\cap\p\A_g$ is also described explicitly. We hope
that this can be helpful for obtaining an estimate of the dimension
of $\tn^h$.

\smallskip
\begin{rem}\label{debarreproof}
O. Debarre has explained to us the following way to easily fix the
proof of the theorem \ref{mainthm} using the results from his work
\cite{debarredecomposes}, without refering to the moduli space of
curves. Such a proof proceeds by induction in $g$. Indeed, just
before ``Quatri\`eme pas'' (p. 701 in \cite{debarredecomposes}), it
is shown that the morphism $S$ is smooth and that $S\to N_g$ is
birational at a general point of $\partial'\tn$. Thus the same must
holds over a general point of $\tn$; moreover, the differential is
injective if and only if the singularity is an ordinary double
point, and thus the result is proven.
\end{rem}

\section{The locus $\tn^{g-1}$}
We will now consider the intersection $\tn\cap N_0'$ and its
eventual relation with $\tn^{g-1}$.

The intersection of the two components of $N_0=\tn+2N_0'$ is studied
in the last section of \cite{debarredecomposes}: it is proven that
their intersection is not reduced for $g\ge 4$, reducible for $g\ge
5$ and irreducible for $g=4$. In a recent paper \cite{genus4} we
proved that in genus $4$ scheme-theoretically
$$
 \tn\cap N_0'=2\tn^3
$$
(for genus 4 the locus $N_0'$ is the Jacobian locus). This was done
by proving an inclusion and checking that the components of the two
varieties have the same degree in the space $\A_4 (4,8)$ (more
precisely, in the projective space $\P^{135}$ containing
$Th(\overline{\A_4(4,8)})$.

It is natural to ask what the situation is for higher $g$. We thus
recall some more notations and results from \cite{debarredecomposes}
and \cite{mu}. Following Mumford, we denote by $\S:=\Sing_{\rm
vert}\T$ the locus of singular points of theta divisors of ppavs.
This is a subvariety of the universal family $$\pi:\X_g\to\A_g.$$
Each component  of $\S$ has codimension $g+1$, cf. \cite{mu} and
\cite{amsp}, and it is locally defined within $\X_g$ by $g+1$
equations
\begin{equation}\label{Sn}
 \t(\tau,z)=0,\qquad \frac{\p \t}{\p z_i}(\t,z)=0\quad \forall i=1,\dots,g.
\end{equation}
Thus set-theoretically $\S$ decomposes into $\Sn$ --- the locus
where $z$ is an even point of order two on the ppav lying on the
theta divisor --- and the remaining component(s) $\S'$. The
following lemma ties in the locus $\tn^{g-1}$ with this description
of the geometry.

\begin{prop}
Set-theoretically we have $\tn^{g-1}=\pi(\Sn\cap\Sing\S)$.
\end{prop}
\begin{proof}
Since the locus $\S$ is given as a subvariety of $\X_g$ by the $g+1$
equations (\ref{Sn}), $\Sing\S$ is the locus where the $g+1$
gradients of these equations, with respect to all the local
coordinates on $\X_g$, i.e. with respect to all $\tau_{ij}$ and
$z_i$, are linearly dependent, i.e. $\Sing\S$ is the locus where the
$(\frac{g(g+1)}{2}+ g)\times (g+1)$ matrix
$$
  \pmatrix
   \frac{\p\t}{\p\tau_{11}}&\ldots&\frac{\p\t}{\p\tau_{gg}}
   &\frac{\p\t}{\p z_1}&\ldots &\frac{\p\t}{\p z_g}\\
   \frac{\p^2\t}{\p z_1\p\tau_{11}}&\ldots&\frac{\p^2\t}{\p z_1\p\tau_{gg}}
   &\frac{\p^2\t}{\p z_1\p z_1}&\ldots &\frac{\p^2\t}{\p z_1\p z_g}\\
   \vdots&\ddots&\vdots&\vdots&\ddots&\vdots\\
   \frac{\p^2\t}{\p z_1\p\tau_{11}}&\ldots&\frac{\p^2\t}{\p z_1\p\tau_{gg}}
   &\frac{\p^2\t}{\p z_g\p z_1}&\ldots &\frac{\p^2\t}{\p z_g\p z_g}\\  \endpmatrix
$$
evaluated at $(\tau,z)$ has rank at most $g$.

If $z$ is an even point of order two, then all $\frac{\p\t}{\p z_i}$
and, by using the heat equation, all $\frac{\p^2\t}{\p
z_i\p\tau_{ij}}$ derivatives at $(\tau,z)$ are equal to zero, and
thus the above matrix becomes
$$
  \pmatrix
   \frac{\p\t}{\p\tau_{11}}&\ldots&\frac{\p\t}{\p\tau_{gg}}&0&\ldots &0\\
   0&\ldots&0&\frac{\p^2\t}{\p z_1\p z_1}&\ldots &\frac{\p^2\t}{\p z_1\p z_g}\\
   \vdots&\ddots&\vdots&\vdots&\ddots&\vdots\\
   0&\ldots&0&\frac{\p^2\t}{\p z_g\p z_1}&\ldots &\frac{\p^2\t}{\p z_g\p z_g}\\
   \endpmatrix.
$$
The rank of this matrix is non-maximal if either all the derivatives
$\frac{\p\t}{\p\tau_{ij}}=0$ for all $1\le i\le j$ (and then by the
heat equation in fact the matrix is completely zero), or if the
$g\times g$ Hessian matrix in the lower right is degenerate, i.e.
exactly if $\tau\in\tn^{g-1}$.
\end{proof}

In the last section of \cite{debarredecomposes} Debarre shows that
$\tn\cap N_0'$ is reducible for $g\ge 5$. In fact he shows that
there is a component $R_g$ of $\tn\cap N_0'$ (the boundary of which
in the partial toroidal compactification is described explicitly)
such that $R_g\subset \pi(\Sn\cap S')$, and that there are other
components $D_g$, which do not lie in $\pi(\Sn\cap S')$.

For a generic ppav in $R_g$ the theta divisor has only one singular
point (of order two) that is the limit of singular points $x$ and
$-x$ of theta divisors of ppav in $N_0'$. For generic points of the
other component(s) $D_g$ of $\tn\cap N_0'$ the theta divisor has a
singular point of order two, and also two other singular points $x$
and $-x$ --- this is exactly to say that $D_g$ does not lie in
$\pi(\Sn\cap S')$.

We now tie in the locus $\tn^{g-1}$ with this picture.
\begin{prop}
The locus $\tn^{g-1}$, as a set, is equal to $\pi(\Sn\cap \S')$.
\end{prop}
\begin{proof}
Indeed, we formally compute (on the level of sets, i.e. with reduced
scheme structure)
$$
\begin{aligned}
 \tn^{g-1}&=\pi(\Sn\cap\Sing\S)=\pi\big(\Sn\cap\Sing(\Sn\cup\S')\big)\\
 &=\pi\Big(\Sn\cap\big((\Sn\cap\S')\cup
 (\Sing\Sn)\cup(\Sing\S')\big)\Big)\\
 &=\pi\big((\Sn\cap\S')\cup\Sing\Sn\big),
\end{aligned}
$$
so what it remains to prove is that
$\pi(\Sing\Sn)\subset\pi(\Sn\cap\S')$. Note, however, that since the
locus $\tn^{g-1}\subset\A_g$ is given locally by two equations ---
some theta constant and its Hessian being zero --- it is purely
codimension two in $\A_g$ (we proved above that it is not
codimension one, since it is not equal to the irreducible divisor
$\tn$, in which it is contained). Therefore it suffices to show that
$\pi(\Sing\Sn)$ has codimension higher than two --- then it cannot
be an irreducible component of $\tn^{g-1}$, and thus must be
contained in $\pi(\Sn\cap S')$. Thus the following lemma finishes
the proof of the proposition.
\end{proof}

\begin{lm}
The codimension of $\pi(\Sing\Sn)$ in $\A_g$ is higher than two.
\end{lm}
\begin{proof}
By definition $\Sn$ is locally given in the universal family $\X_g$
by the following $g+1$ equations:
$$  \t(\tau,z)=0,\qquad  z=(\tau\e+\de)/2$$
where $[\e,\de]$ is an even theta characteristic. Thus the locus
$\Sing\Sn$ is where the $g+1$ gradient vectors (with respect to all
local coordinates on $\X_g$, i.e. with respect to both $\tau$ and
$z$) of these defining equations are linearly dependent. These
gradients form the matrix
$$
 \pmatrix
  \frac{\p\t}{\p\tau_{11}}&\ldots&\frac{\p\t}{\p\tau_{gg}}
   &\frac{\p\t}{\p z_1}&\ldots &\frac{\p\t}{\p z_g}\\
  -\e_1/2&\ldots&0&1&\ldots&0\\
  \vdots&\ddots&\vdots&\vdots&\ddots&\vdots\\
  0&\ldots&-\e_g/2&0&\ldots&1
 \endpmatrix
$$
When we evaluate the derivatives of the theta function at the even
point $(\tau\e+\de)/2$ of order two, the derivatives $\frac{\p\t}{\p
z_i}$ all vanish. Thus the matrix only has non-maximal rank if all
the derivatives $\frac{\p\t}{\p\tau_{ij}}$ are zero for all $1\le
i\le j\le g$, i.e. if the corresponding theta constants $\tt\e\de$
vanishes at $\tau$ to order at least two. A naive dimension count
would predict a high codimension for such a condition; here is an
easy observation showing that this codimension is greater than two.

Indeed, by the heat equation the theta constants vanishes to order
at least two only if the theta function $\tt\e\de(\tau,z)$ vanishes
at $z=0$ to order at least 4 (the third $z$-derivatives at zero are
all zero by parity). In \cite{te} Teixidor i Bigas studies the locus
of curves having such a theta characteristic --- this is $\M_g^3$ in
her notations
--- and shows that ${\rm codim}_{\M_g}\M_g^3>2$ (her results are
actually better than this, but this is all we need). But if a locus
in $\A_g$ has non-empty intersection with $\M_g\subset\A_g$, and the
codimension of its intersection with $\M_g$ within $\M_g$ is at
least $n$, then the codimension of the locus itself in $\A_g$ is at
least $n$, cf. \cite{cpv}. Rigorously we can only use this for
smooth varieties, while $\A_g$ is an orbifold, but we can pass to
the finite covering $\A_g (4,8)$ that is smooth. Alternatively, we
can use the result from \cite{amsp} that if the theta divisor has a
point of multiplicity greater than two, then the period matrix
belongs to the Andreotti-Mayer locus $N_1$ that has codimension at
least 3 in $\A_g$.
\end{proof}
\begin{cor} As sets (i.e. with reduced scheme structure)
we have the inclusion $\tn^{g-1}\subset(\tn\cap N_0')$.
\end{cor}
\begin{proof}
Indeed, by definition we have $\tn=\pi(\Sn)$, and $N_0'=\pi(\S')$,
and thus
$\tn^{g-1}=\pi(\Sn\cap\S')\subset(\pi(\Sn)\cap\pi(\S'))=\tn\cap
N_0'$.
\end{proof}
Furthermore, we can now describe more precisely the locus $R_g$
introduced by Debarre.
\begin{prop}
We have the equality of sets $R_g=\tn^{g-1}=\pi(\Sn\cap\S')$.
\end{prop}
\begin{proof}
Recall that $R_g$ is defined in \cite{debarredecomposes} by first
studying the boundary of the locus $\tn\cap N_0'$ in the partial
toroidal compactification $\overline{\A_g}^1$, choosing a certain
explicitly defined component $\p R_g$ of this boundary, and then
arguing that $\p R_g$ must be the boundary of some locus
$R_g\subset\A_g$. It then follows that $R_g\subset\pi(\Sn\cap\S')$,
and it is shown in \cite{debarredecomposes} that no other component
of $\tn\cap N_0'$ is contained in $\pi(\Sn\cap\S')$.

Since we know the equality $\tn^{g-1}=\pi(\Sn\cap\S')$, it follows
that $\pi(\Sn\cap\S')$ is purely codimension two, and since it is
contained in $\tn\cap N_0'$, which is also purely codimension two,
$\pi(\Sn\cap\S')$ must be the union of a number of irreducible
components of $\tn\cap N_0'$. But then since this union cannot
contain any of the components of $D_g$, we must have
$\pi(\Sn\cap\S')=R_g$.
\end{proof}
The equality proved above is set-theoretic. Since $R_g$ and
$\tn^{g-1}$ have the same intersection with the boundary $\partial
\A_g$, we also get
\begin{cor}
Up to embedded subvarieties $\tn^{g-1}=R_g$ as schemes.
\end{cor}

\section*{Acknowledgements}
We would like to thank Ciro Ciliberto and Gerard van der Geer for
discussions regarding the relations between $\tn^{g-1}$ and $\tn\cap
N_0'$, and Olivier Debarre for reading a preliminary version of this
note and for comments on the relation of $\tn^{g-1}$ and $\tn$ (see
the introduction). Special thanks to Marco Manetti for suggesting an
elegant proof of lemma \ref{mainlm}.

\end{document}